\documentclass[leqno,12pt]{article}

\usepackage{latexsym}
 \usepackage{graphicx}
\usepackage{amsmath}
\usepackage{amssymb}
\usepackage{algorithm}
\usepackage{algorithmic}
\usepackage{cite}
\usepackage{color}

\usepackage[margin=1.3in, top=1.3in, bottom=1.3in]{geometry}

\newcommand{\nc}{\newcommand}
\nc{\nt}{\newtheorem}
\nt{thm}{Theorem}[section]
\nt{cor}[thm]{Corollary}
\nt{prop}[thm]{Proposition}
\nt{lem}[thm]{Lemma}
\nt{defn}[thm]{Definition}
\nt{rem}[thm]{Remark}
\nt{exa}[thm]{Example}
\nt{ass}[thm]{Assumption}
\nt{alg}[thm]{Algorithm}
\nt{que}[thm]{Question}
\nt{con}[thm]{Conjecture}
\nc{\ip}[2]{\mbox{$\langle #1,#2 \rangle$}}
\nc{\pf}{\noindent{\bf Proof\ \ }}
\nc{\finpf}{\hfill{$\Box$}\linespace}
\nc{\linespace}{\vspace{\baselineskip} \noindent}
\nc{\R}{{\bf R}}
\nc{\Rn}{{\bf R}^n}
\nc{\Sn}{{\bf S}^n}
\nc{\lm}{\lambda_{\mbox{\rm\scriptsize min}}}
\nc{\bx}{\bar{x}}
\nc{\e}{\epsilon}
\nc{\cl}{\mbox{\rm cl}\,}
\nc{\conv}{\mbox{\rm conv}\,}

\def\tto{\;{\lower 1pt \hbox{$\rightarrow$}}\kern -12pt
           \hbox{\raise 2.8pt \hbox{$\rightarrow$}}\;}
\newenvironment{myequation}{\setcounter{equation}{\value{thm}}
   \begin{equation}}{\addtocounter{thm}{1}\end{equation}}
\newenvironment{myeqnarray}{\setcounter{equation}{\value{thm}}
    \begin{eqnarray}}{\setcounter{thm}{\value{equation}}\end{eqnarray}}

\nc{\bmye}{\begin{myequation}}
\nc{\emye}{\end{myequation}}

\begin{document}
\title{
Rescaling nonsmooth optimization using BFGS and Shor updates
}
\author{
 Jiayi Guo
\thanks{ORIE, Cornell University, Ithaca, NY 14853, U.S.A.
\texttt{jg826@cornell.edu}.}
\and 
A.S. Lewis
\thanks{ORIE, Cornell University, Ithaca, NY 14853, U.S.A.
\texttt{people.orie.cornell.edu/aslewis}.
Research supported in part by National Science Foundation Grant DMS-1613996.}
}
\maketitle

\begin{abstract}
The BFGS quasi-Newton methodology, popular for smooth minimization, has also proved surprisingly effective  in nonsmooth optimization.  Through a variety of simple examples and computational experiments, we explore how the BFGS matrix update improves the local metric associated with a convex function even in the absence of smoothness and without using a line search.  We compare the behavior of the BFGS and Shor r-algorithm updates.
\end{abstract}
\medskip

\medskip

\medskip

\section{Introduction}
We consider unconstrained minimization methods for a function $f \colon \Rn \to \R$.  Our aim is to explore basic theory, so for simplicity we assume throughout that $f$ is convex and everywhere finite, even though many of the algorithms we consider are also interesting for functions that may be nonconvex or extended-valued.

Since the 1970's, an extensive literature has documented the powerful properties of the BFGS (Broyden-Fletcher-Goldfarb-Shanno) update in quasi-Newton minimization algorithms.  The update accumulates information about the curvature of the objective, allowing, like Newton's method, a beneficial transformation of the space.  Remarkably, this process seems to help reliably even when the objective is nonsmooth \cite{BFGS}.  In this work  we try to illuminate this phenomenon.

When studying the BFGS algorithm in the context of nonsmooth optimization, an interesting point of comparison is the Shor r-algorithm \cite{shor}.  Shor's method also uses a quasi-Newton-like transformation, simpler than BFGS, but without satisfying the secant condition standard in smooth optimization \cite{nocedal_wright}.  The algorithm is challenging to analyze \cite{shor-speed}, and hard to implement systematically in practice, although there have been promising attempts \cite{kappel-kuntsevich}.  One fundamental difficulty is how to incorporate into the method a systematic line search.  By contrast, a simple line search (satisfying the standard weak Wolfe conditions) is easy to incorporate into a nonsmooth BFGS algorithm, and seems very successful in practice \cite{BFGS}.

Unfortunately, we lack almost any theoretical insight into the benefits of the BFGS (or Shor) update in nonsmooth optimization.  The interplay of a line search with quasi-Newton updates complicates the question still further.  In this work, we therefore try to isolate the behavior of the BFGS update, in particular, and try to understand its beneficial effects {\em with no line search\/}.

On simple random nonsmooth convex optimization problems, a BFGS method can dispense almost entirely with the usual line search and seemingly still reliably succeed.  More precisely, if the objective $f$ is smooth at the current iterate $x \in \Rn$ (as holds generically) and the current BFGS matrix is $H$ ($n$-by-$n$ and positive-definite), then traditionally we calculate $x_+ = x - tH\nabla f(x)$, where the stepsize $t > 0$ is chosen through the line search, apply a standard BFGS update formula to $H$, and update $x = x_+$.  However, a more rudimentary idea is simply to choose $t=1$, and update
$H$, but only update $x$ if the step generates descent:  $f(x_+) < f(x)$.  We refer to this stripped-down method as {\em linesearch-free BFGS}.

\begin{figure} 
\begin{center}
\includegraphics[width=80mm]{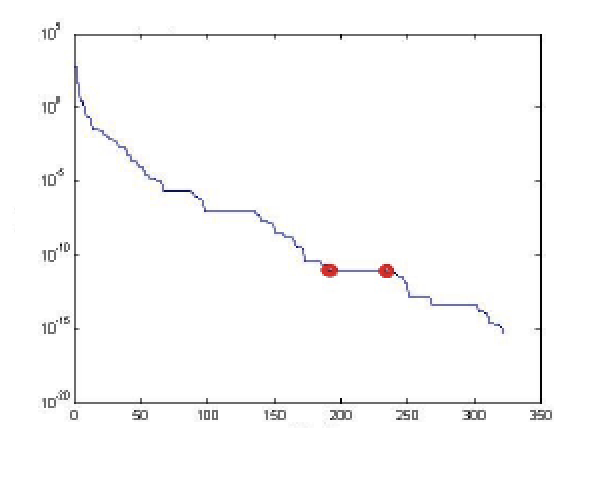}
\end{center}
\caption{A typical run of the linesearch-free BFGS method for a nonsmooth $f$ on $\R^5$, plotting the value $f(x_k) - \min f$ against the iteration count $k$.  Flat segments (such as between the red dots) indicate BFGS updating without iterate updating.}
\label{four-random}
\end{figure}

In Figure \ref{four-random}, we illustrate the typical performance of the linesearch-free BFGS method on a simple example, where the objective function $f \colon \R^5 \to \R$ is the maximum of four random strictly convex quadratics.  Somewhat surprisingly, the method minimizes $f$ reliably, generating a sequence of iterates that appear to converge linearly to the minimizer.  Experiments such as these prompt our quest for insight into the BFGS update for nonsmooth functions.

\section{Space dilation via Shor and BFGS}
\label{space-dilation}
We begin our exploration with a simpler method: the Shor r-algorithm for minimizing the function $f$.  At each iteration we consider the current iterate $x \in \Rn$, and a current subgradient $g \in \partial f(x)$.  The classical subgradient method takes a step from $x$ in the direction $-g$.  Shor \cite[Section 3.6]{shor} proposed accelerating this idea by successively rescaling the space using a current $n$-by-$n$ matrix $V$, initially equal to the identity matrix $I$.  We update the iterate 
via
\bmye \label{shor-factors}
s = -V^T V g; ~~ x_+ = x + ts;
\emye
the stepsize $t > 0$ being chosen through some kind of line search.  At the new iterate $x_+$ we then find a new subgradient $g_+ \in \partial f(x_+)$, define a unit vector $e \in \Rn$ by
\[
e = V(g - g_+); ~~ e = \frac{e}{\|e\|};
\]
update the matrix via
\[
W = I  - \frac{e e^T}{2\|e\|^2};  ~~ V_+ = WV;
\]
update $x=x_+$; $g=g_+$; $V=V_+$;
and repeat.  (The factor ``2''  that appears in the denominator in the definition of $W$ has no special significance and could be replaced by any constant greater than $1$.)  Shor described his method as one of ``space dilation'':  after making a current change of variables $x=V^T y$, the unit vector $e$ lies in the direction of the difference of two successive subgradients of the objective function $y \mapsto f(V^T y)$, and the transformation $W$ dilates the space in this direction.

Consider the canonical special case of minimizing a sublinear function:
\[
f(x) = \max_{h \in Q} h^T x,
\]
for a nonempty compact set of nonzero vectors $Q \subset \Rn$.  Then $x=0$ is nonoptimal if and only if there exists a descent direction:  a vector $d \in \Rn$ and a scalar $\alpha < 0$ such that $h^T d \le \alpha$ for all $h$ in $Q$.   This condition states that a hyperplane normal to $d$ separates zero from $Q$ (or equivalently its convex hull $\mbox{conv}\,Q$).

We could apply Shor's method, seeking to minimize the function $f$ starting (and remaining) at the point zero, and terminating once we find a descent direction.  More precisely, at each iteration the current iterate is $x=0$ and the current subgradient $g$ lies in the set $Q$.  We choose the stepsize $t=1$, terminate if $f(s) < 0$, and otherwise choose a new subgradient $g_+ \in Q$ to maximize the inner product $s^T g_+$. (As we discuss and motivate in Section \ref{nonsmooth-updating}, this choice of $g_+$ correctly models the function in the search direction: $f(s) = s^T g_+$.)  We then update the matrix $V$ and the subgradient $g$, maintain $x=0$, and repeat.

Following the change of variables we introduced above, if we define $h=Vg$ and $p=Vg_+$, we arrive at the following simple procedure for separating a set $Q$ from zero, relying only on a linear optimization oracle over $Q$.

\begin{alg}[Shor update method for $0 \in \mbox{conv}\,Q$]
\label{shor-update-method}
{\rm
\begin{algorithmic}
\STATE
\STATE  Choose $h \in Q$; $V = I$;
\WHILE{not done}
\STATE  find a minimizer $p$ of $\ip{\cdot}{h}$ over $Q$;
\IF{$p^Th > 0$}
\STATE  terminate with $V^T h$ ``normal to separating hyperplane'';
\ENDIF
\STATE  $e = h-p$; $W = I - \frac{ee^T}{2\|e\|^2}$; $Q = WQ$; $V = WV$; $h = Wp$;
\ENDWHILE
\end{algorithmic}
}
\end{alg}

\noindent
Geometrically, the procedure tests whether the current vector $h \in Q$ is normal to a hyperplane separating $Q$ from zero, and if not applies to the set $Q$ a simple linear transformation $W$, a symmetric rank-one perturbation of the identity.

In its intuitive simplicity and apparent versatility, the procedure has a certain appeal.  Furthermore, experiments on small examples suggest it typically works:  when zero lies outside the convex hull of $Q$, the procedure terminates, and otherwise the vector $h$ converges to zero.  Especially simple is the case when $Q$ is a finite set of nonzero vectors, in which case we seek to separate the polytope $\mbox{conv}\,Q$ from zero, or equivalently, find a solution $d$ for the homogeneous system of inequalities $h^Td<0$ for all $h \in Q$ (a core problem of linear programming).  The oracle --- linear optimization over $Q$ --- is then particularly easy.

Each iteration is computationally simple, involving just elementary operations on all the vectors in $Q$.  The procedure is not as ``elementary'' as methods like the Perceptron Algorithm, that in particular preserve sparsity, being closer in spirit to rescaled perceptron methods \cite{beloni-freund-vempala}.  As we discuss later, it also has some formal similarities with versions of the Ellipsoid Algorithm.  

\begin{figure} 
\begin{center}
\includegraphics[width=100mm]{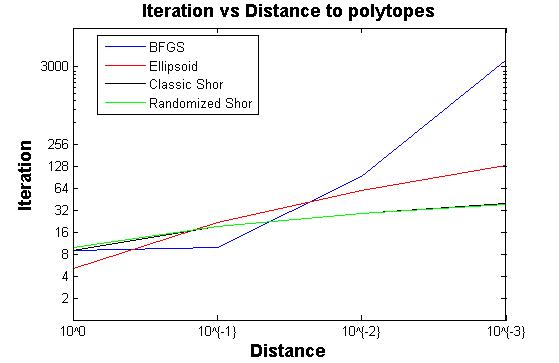}
\end{center}
\caption{Separating a point from a polytope:  mean of number of required iterations to terminate from a random start.}
\label{polytope-separation}
\end{figure}

\begin{figure} 
\parbox{75mm}
{
\begin{center}
\includegraphics[width=75mm]{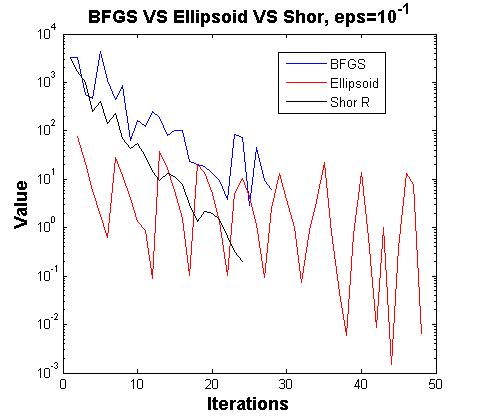}
\end{center}
}
\hfill
\parbox{75mm}
{
\begin{center}
\includegraphics[width=75mm]{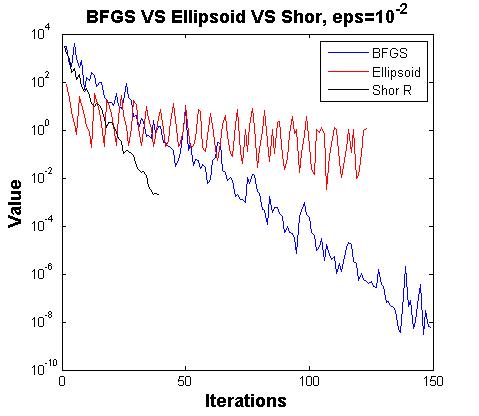}
\end{center}
}
\caption{Separating a point from a polytope:  typical trajectories.}
\label{polytope-separation2}
\end{figure}

To illustrate, consider the following example in dimension $n=5$.  For each index $j$, denote the corresponding unit vector in $\R^5$ by $e_j$.   Define vectors $a_j = 4^j e_j$, along with a convex combination $p = (\sum 4^{-j})^{-1} \sum e_j$.  Fix a parameter $\epsilon > 0$, and let the set $Q$ consist of the points $a_j - (1+\epsilon)p$ (for each $j$) along with $-p$.  Geometrically, $Q$ consists of the vertices of an irregular simplex.  When $\epsilon$ is small the point zero is outside the simplex but close to one of the facets, making the problem ill-posed.  Figure \ref{polytope-separation} plots the number of iterations needed by several algorithms to find a separating hyperplane, averaged over the starting point in $Q$, as a function of the parameter $\epsilon$.  The figure compares this Shor update method (labeled ``Classic Shor'') with several other algorithms discussed below:  a randomized Shor method, a BFGS method, and a version of the ellipsoid algorithm.   As we see from the plot, on this small and simple example, Shor updating is reliable, terminating after a couple of dozen iterations even with $\epsilon = 10^{-3}$.  Figure \ref{polytope-separation2} shows some typical trajectories.

A simple nonpolyhedral problem seeks to separate a point $c$ from an ellipsoid in $\Rn$.   If we describe the ellipsoid as $AB$, where $B \subset \Rn$ is the closed unit ball and $A$ is an invertible $n$-by-$n$ matrix, then we seek a normal vector $z \in \Rn$ to a separating hyperplane, or in other words a solution of the inequality $\|A^T z\| < c^T z$.  If $Q$ is the boundary of the ellipsoid $AB-c$, we arrive at following simple procedure (involving no matrix inversion):  if it terminates, the output vector $z$ solves our problem.
\newpage

\begin{alg}[Shor updating to separate point $c$ from ellipsoid $AB$] 
\label{simple}
{\rm
\begin{algorithmic}
\STATE
\STATE  Choose unit $x \in \Rn$; $V=I$;
\WHILE{not done}
\STATE  $h = Ax-c$; $y = -A^T h$; $y = \frac{y}{\|y\|}$; $p = Ay-c$;
\IF{$p^T h > 0$}
\STATE  terminate with $z=V^T h$ ``normal to separating hyperplane'';
\ENDIF
\STATE  $e = h-p$; $W = I - \frac{ee^T}{2\|e\|^2}$; $A = WA$; $V = WV$; $c = Wc$; $x = y$;
\ENDWHILE
\end{algorithmic}
}
\end{alg}

\noindent
(In the notation of Algorithm \ref{shor-update-method}, the current iterate $h$ is $Ax-c$ for some unit vector $x$, and we compute $p$ by minimizing $\ip{p}{h}$ with $p = Ay-c$, over unit vectors $y$.)

\begin{figure} 
\parbox{75mm}
{
\begin{center}
\includegraphics[width=75mm]{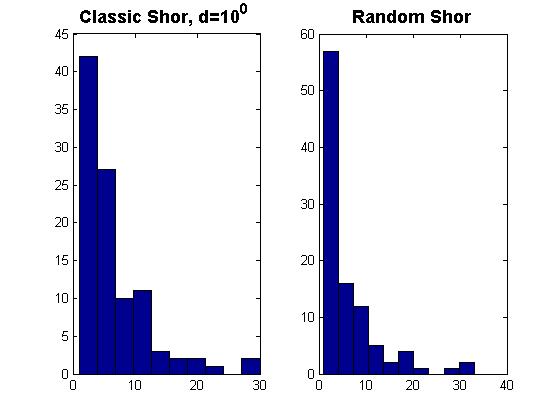}
\end{center}
}
\hfill
\parbox{75mm}
{
\begin{center}
\includegraphics[width=75mm]{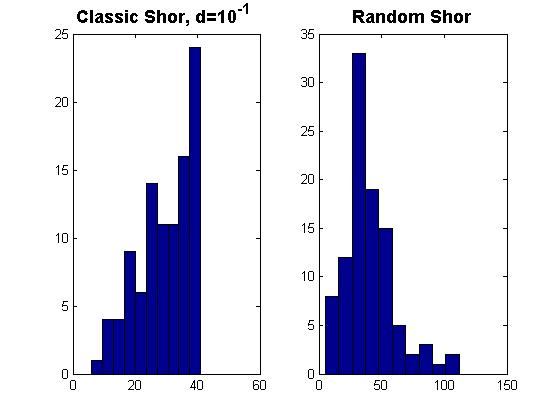}
\end{center}
}
\caption{Shor updating to separate a point from an ellipsoid in $\R^5$.  Histograms of number of required iterations to terminate, for 100 random examples.}
\label{shor-tests}
\end{figure}

We illustrate this idea in dimension $n=5$, for a diagonal matrix $A$ with diagonal 
$[1 ~ 10 ~ 10^2 ~ 10^3 ~ 10^4]$.  We generate a hundred instances by choosing the vector $c= (1+d)Au$, where $u \in \R^5$ is a random unit vector, and we set the scalar $d$ (which controls the ill-posed of the instance) to be both $1$ and $10^{-1}$ to illustrate the effect of ill-posedness.  The figures below (Figure \ref{shor-tests}) plot histograms of the number of instances requiring certain numbers of  iterations to terminate.

As we see from the plots, on these random examples on a low-dimensional ellipsoid, moderately ill-conditioned and well separated from zero, this Shor updating method works reliably.  It typically finds a separating hyperplane after a couple of dozen iterations.  Not surprisingly, the required number of iterations grows as the parameter $d$ (and hence distance to ill-posedness) shrinks:  examples with $d=10^{-2}$ may need more than $100$ iterations, and $d=10^{-3}$ may need more than $1000$.  The method remains viable as the dimension grows.  When the matrix $A$ has diagonal entries 
$[1 ~ 10 ~ 10^2 ~ \cdots ~ 10^9]$, with $d=10^{-1}$, a hundred random instances all terminated in less than 200 iterations.

\begin{figure} 
\begin{center}
\includegraphics[width=100mm]{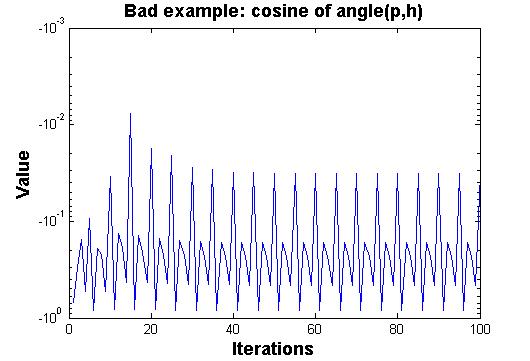}
\end{center}
\caption{Cyclic behavior of the angle between $p$ and $h$ during Shor updating.}
\label{bad-example}
\end{figure}

Bolstered by such random experiments, where Shor updating systematically succeeds, we might hope for a simple proof validating Algorithm \ref{shor-update-method}, and thereby some insight into the Shor r-algorithm.  Sadly, while the procedure is simple, its behavior is not:  sporadically, {\em it can fail}.  For example, numerical experiments with the ellipsoid separation procedure (Algorithm \ref{simple}) revealed that on the small example in $\R^2$ defined by
\[
A = 
\left[
\begin{array}{cc}
1 & 0 \\
0 & 10 
\end{array}
\right],
~~
v = -
\left[
\begin{array}{c}
10 \\
39 
\end{array}
\right],
~~
c = (1+10^{-2})A \frac{v}{\|v\|},
\]
a thousand iterations do not suffice for termination.  Furthermore, the failure seems unambiguous:  after a few steps, iterations seem to behave cyclically, with a period of five iterations. In particular, the cosine of the angle between the vectors $p$ and $h$ is bounded above by $-1/100$, so the termination criterion always fails.  
Figure \ref{bad-example} plots this cosine for the first hundred iterations.

Sporadic failures notwithstanding, the numerical evidence suggests that the Shor update in Algorithm \ref{shor-update-method} improves the geometry of the problem in some average sense, perhaps analogous to randomized algorithms for convex programming like \cite{dunagan-vempala,beloni-freund-vempala}.  To try to isolate this effect, we consider a simple modification of Algorithm \ref{shor-update-method}.  The original method remembers the new element $p \in Q$ at the end of an iteration, and uses it, transformed as $h=Wp$, to start the next iteration.  The modified method below forgets $p$ after the iteration, starting each new iteration afresh by simply choosing $h$ to optimize over $Q$ in a random direction.
\newpage

\begin{alg}[Randomized Shor for $0 \in \mbox{conv}\,Q$]
{\rm
\begin{algorithmic}
\STATE
\STATE  $V = I$;
\WHILE{not done}
\STATE	choose random $u \in \Rn$;
\STATE	find a minimizer $h$ of $\ip{\cdot}{u}$ over $Q$;
\STATE  find a minimizer $p$ of $\ip{\cdot}{h}$ over $Q$;
\IF{$p^Th > 0$}
\STATE  terminate with $V^T p$ ``normal to separating hyperplane'';
\ENDIF
\STATE  $e = h-p$; $W = I - \frac{ee^T}{2\|e\|^2}$; $Q = WQ$; $V = WV$;
\ENDWHILE
\end{algorithmic}
}
\end{alg}

\noindent
We could, for example, distribute the random vector $u$ normally.  In the special case of the ellipsoid separation procedure, we arrive at a Randomized Algorithm \ref{simple}, where $x$, rather than equalling $y$, is just the normalized vector $A^T u$ for a random vector $u$.  In the figure, we compare the results for the randomized procedure with those for the original ``classic'' procedure, and like that procedure, it seems reliable on small random examples.  On the one hand, we observe no failures.  On the other hand, shrinking the ill-posedness parameter $d$ seems to slow the randomized procedure more than the original version.  With $d=10^{-2}$ (not shown in the figure), the original version terminates in every instance within around a hundred iterations, whereas the randomized version often takes thousands.  In summary, reusing the previous element of the set $Q$ at each iteration seems to accelerate the procedure.

We noted in the introduction that the BFGS method, as a general-purpose nonsmooth optimization tool, shows more promise than the Shor r-algorithm, and is quite successful in practice \cite{BFGS}.
A single modification to each of the space-dilation algorithms above transforms the motivation from the Shor update to the BFGS quasi-Newton update.  Specifically, as we explain in Section \ref{cholesky-factors}, we simply change the updating transformation from
\[
W ~=~ I - \frac{ee^T}{2\|e\|^2}
\]
to 
\[
W ~=~ I - \Big(\frac{e}{h^Te} - \frac{h}{\|h\|\sqrt{h^Te}}\Big) h^T
\]
A geometric interpretation of the resulting BFGS-based procedure is almost identical to that for the Shor updating procedure.  The only difference is that the matrix $W$ transforming the space, while still a rank-one perturbation of the identity, is now no longer symmetric.  As with the Shor update, any convergence theory for this BFGS procedure seems elusive, but its simplicity and apparent effectiveness are intriguing.

\begin{figure} 
\begin{center}
\includegraphics[width=100mm]{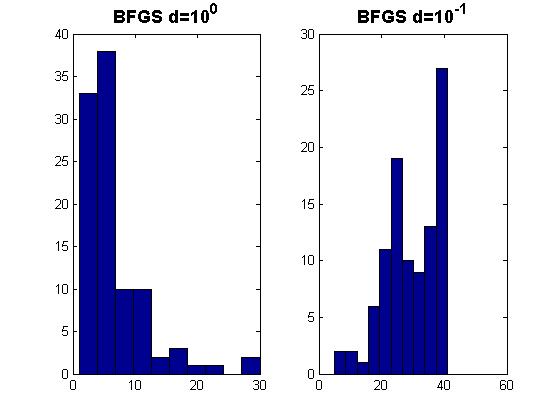}
\end{center}
\caption{BFGS updating to separate a point from an ellipsoid in $\R^5$.  Histograms of number of required iterations to terminate, for 100 random examples.}
\label{bfgs-ellipsoid}
\end{figure}

On polyhedral separation examples (Figure \ref{polytope-separation}), the BFGS method is successful but seems slower than the Shor method as the ill-posedness parameter $\epsilon$ shrinks.  The previous ellipsoid separation examples suggest a similar comparison (Figure \ref{bfgs-ellipsoid}):  the Shor method seems faster as the ill-posedness parameter shrinks (although the randomized version seems slower), and sporadic failure is a possibility.

\section{The BFGS update}
To begin a more careful discussion of the BFGS update, we first recall the classical idea of Newton's method as a steepest descent method in a local metric. 
Suppose the function $f$ is ${\mathcal C}^2$-smooth, and consider a point $x \in \Rn$ at which the Hessian $\nabla^2 f(x)$ is positive definite.  If we denote the gradient $\nabla f(x)$ by $g$, then
the unit steepest descent step, with respect to the  Euclidean norm, is the minimizer of the linear approximation $g^T s$ over the unit ball 
$\{s : \|s\| \le 1 \}$, namely $s = -\frac{1}{\|g\|}g$ (assuming $g \ne 0$).  If instead we minimize over the ball 
$\{s : s^T \nabla^2 f(x) s \le 1\}$ corresponding to a natural local metric associated with $f$ at $x$, we instead arrive at the Newton step $s = -\nabla^2 f(x)^{-1}g$.  In stark contrast to the steepest descent step, taking the Newton step from $x$ and iterating --- Newton's method --- rapidly reduces the objective value, at least close to a minimizer $\bar x$ with $\nabla^2 f(\bar x)$ positive definite.  If the initial point is far from $\bar x$, a simple backtracking line search along the direction of the Newton step can guarantee progress into the neighborhood where the unit Newton step is acceptable.

Turning to quasi-Newton methods, the classical idea is to use an approximation $H$ in the place of the inverse Hessian $\nabla^2 f(x)^{-1}$, updated after each step.  In particular, the BFGS (Broyden-Fletcher-Goldfarb-Shanno) update, assuming, like Newton's method, a unit step, uses a matrix $H$ in the set 
$\Sn_{++}$ of positive-definite $n$-by-$n$ symmetric matrices.  We update $H$ as follows:
\begin{myeqnarray}
s = -Hg, ~~ x_+ = x+s, & &  g_+ = \nabla f(x_+), ~~ y = g_+ - g  \label{update1} \\
V = I - \frac{sy^T}{s^Ty}, & & H_+ = VHV^T + \frac{ss^T}{s^Ty}. \label{update2}
\end{myeqnarray}
In this update, we assume that the quantity $s^Ty$ is strictly positive:  it must be nonnegative, by convexity, but unless $f$ is strictly convex, it may be zero.  For future reference, we make the following definition.

\begin{defn} \label{unit-step}
{\rm 
Given a ${\mathcal C}^1$-smooth convex function $f \colon \Rn \to \R$ and a point $x$ in $\Rn$, denote the gradient $\nabla f(x)$ by $g$.  The {\em unit-step BFGS update} is the map $\mbox{BFGS}_{f,x} \colon \Sn_{++} \to \Sn_{++}$ defined by 
$\mbox{BFGS}_{f,x}(H) = H_+$ for a matrix $H \in \Sn_{++}$, where the matrix $H_+$ is given by equations (\ref{update1}) and (\ref{update2}).  If $s^Ty = 0$ in equation (\ref{update1}), the matrix $H_+$ is undefined.
}
\end{defn}

Like Newton's method, if the initial point $x$ and approximation $H$ are close to the minimizer $\bar x$ and inverse Hessian $\nabla^2 f(\bar x)^{-1}$ respectively, then updating $x = x_+$, $g = g_+$, 
$H = H_+$,  and iterating, rapidly reduces the objective value \cite[Thm 8.6]{dennis1977quasi}.  Far from $\bar x$, a line search again can guarantee progress into the neighborhood where the unit step is acceptable, resulting in an algorithm with good global {\em and} local convergence properties:  see \cite{nocedal_wright} for more details and an extended discussion of the enduringly popular BFGS method.  Again like Newton's method, rather than thinking of $H$ as an inverse Hessian approximation, we can instead associate it with a local metric at $x$, and thereby interpret the BFGS as a variable metric method.  This viewpoint better suits our current development, where the objective function $f$ may not be smooth.

For {\em nonsmooth} optimization, although supported by little theory, extensive computational experiments suggest that the BFGS method can also be surprisingly effective \cite{BFGS}.  Under reasonable conditions, with a suitably randomized initial point, the function $f$ is smooth at all points encountered by the method, so the update equations (\ref{update1}) and (\ref{update2}) make sense.  In general, as in the smooth case, it is crucial to incorporate a suitable line search, scaling the step $s$ defined in equation (\ref{update1}) at the outset of the update.  In particular, the experiments in \cite{BFGS} rely on a weak Wolfe line search, ensuring both a sufficient decrease condition on the new objective value $f(x_+)$ and a curvature condition.  Classically, these conditions serve multiple purposes.  The sufficient decrease condition is important in convergence proofs (although in practice simply ensuring decrease typically seems to suffice).  The curvature condition guarantees in particular the condition $s^Ty > 0$, which in turn ensures that the update $H_+$ is positive definite (although for strictly convex objectives $f$ this property is automatic).  Well-known self-correcting properties of the BFGS update, allowing recovery from badly scaled approximations $H$, are thought to depend heavily on the line search \cite{nocedal_wright}.

Important as it is, the line search complicates the already challenging task of understanding how the BFGS update improves the local metric, especially in the nonsmooth case.  The line search may sometimes be irrelevant, as occurs asymptotically in the smooth case, for example.  We therefore ask:  
\begin{itemize}
\item
What can we learn simply from the unit-step BFGS update (\ref{update1}) and (\ref{update2}), {\em with no line search}?   
\item
In particular, when does repeated application of the unit-step BFGS update at a fixed point generate a descent step?  
\item
Can unit-step BFGS updating underly a nonsmooth minimization algorithm?
\end{itemize}

\section{BFGS updating for smooth functions}
Although our primary interest is in nonsmooth functions, we begin with the simplest smooth case.
Consider a ${\mathcal C}^1$-smooth strictly convex function $f \colon \Rn \to \R$, a point $x \in \Rn$ at which the gradient $g = \nabla f(x)$ is nonzero, and an initial matrix $H \in \Sn_{++}$.  Consider the following procedure, which repeatedly applies the unit-step BFGS update (Definition \ref{unit-step}) at the fixed point $x$ until it generates a descent step.
\newpage

\begin{alg}[BFGS updating for smooth function $f$] 
\label{descent}
{\rm
\begin{algorithmic}
\STATE
\WHILE{$f(x - Hg) \ge f(x)$}
\STATE  $H = \mbox{BFGS}_{f,x}(H)$;
\ENDWHILE
\end{algorithmic}
}
\end{alg}

\noindent
(Strict convexity of $f$ ensures that the update is always well-defined.)  When must this iteration terminate?  We begin with the one-dimensional case.

\begin{thm}
BFGS updating (Algorithm \ref{descent}) terminates for any ${\mathcal C}^1$-smooth strictly convex function $f \colon \R \to \R$ at any noncritical point $x$.
\end{thm}

\pf
We argue by contradiction.  Suppose without loss of generality $g = f'(x) = 1$ and the iteration does not terminate.  The matrix $H$ is now simply a scalar $h>0$, and we have $s = -h$, $x_+ = x - h$, 
$g_+ = f'(x-h)$, $y = f'(x-h) - 1$, and $V=0$.  We then update:
\[
h \leftarrow h_+ = \frac{h}{1 - f'(x-h)}.
\]
By assumption, $f(x-h) \ge f(x)$, so $f'(x-h) < 0$.  Hence $h_+ < h$ at every iteration, so $h$ decreases to some limit $\bar h \ge 0$.  By continuity we have $f'(x-\bar h) \le 0$.  If $f'(x-\bar h) = 0$, then we obtain a contradiction, since then the points $x-h$ approach the minimizer $x - \bar h$ so eventually 
$f(x-h) < f(x)$.  Hence in fact we have $f'(x-\bar h) < 0$.  But now we obtain the contradiction
\[
h_+ = \frac{h}{1 - f'(x-h)} \to \frac{\bar h}{1 - f'(x-\bar h)} < \bar h.
\]
This completes the proof.
\finpf

Using standard theory from the classical quasi-Newton literature \cite{dennis1977quasi,griewank}, we next work towards an analogous result for multivariate quadratic functions.  Consider first the special case
$f(x) = \frac{1}{2}\|x\|^2$ (for $x \in \Rn$).  In that case a quick calculation shows that Algorithm \ref{descent} becomes the following.

\begin{alg}[BFGS updating for $\frac{1}{2}\|\cdot\|^2$] 
\label{bfgs-norm}
{\rm
\begin{algorithmic}
\STATE
\WHILE{$\|x - Hx\| \ge \|x\|$}
\STATE	$s = -Hx$; $z = \frac{s}{\|s\|}$; $H = (I-zz^T)H(I-zz^T) + zz^T$;
\ENDWHILE
\end{algorithmic}
}
\end{alg}

In fact the case of a general strictly convex quadratic function $f(x) = \frac{1}{2}\|Rx\|^2$ (for $x \in \Rn$), where the $n$-by-$n$ matrix $R$ is invertible, follows immediately from this special case.  The change of variables $\hat x = Rx$ and $\hat H = RHR^T$ shows, after some algebra, that the BFGS updating algorithm is essentially identical to Algorithm \ref{bfgs-norm}:

{\rm
\begin{algorithmic}
\STATE
\WHILE{$\|\hat x - \hat H \hat x\| \ge \|\hat x\|$}
\STATE	$s = -\hat H\hat x$; $z = \frac{s}{\|s\|}$; $\hat H = (I-zz^T)\hat H(I-zz^T) + zz^T$;
\ENDWHILE
\end{algorithmic}
}
\smallskip

\noindent
To proceed, we begin with some geometry in the Euclidean space $\Sn$ of $n$-by-$n$ symmetric matrices with the inner product defined by $\ip{X}{Y} = \mbox{trace}(XY)$, for matrices $X,Y \in \Sn$.  We start with a tool whose proof is immediate.

\begin{lem}
\label{orthogonal}
Consider any matrix $P \in \Sn$ satisfying $P^2 = P$.  For any matrix $X \in \Sn$, the matrix 
$X_+ = PXP$ is orthogonal to the matrix $X_+ - X$.
\end{lem}

\noindent
For the next step, we denote the smallest eigenvalue of a matrix $H \in \Sn$ by $\lm(H)$.

\begin{lem} \label{tool1}
For any unit vector $z \in \Rn$, and any matrix $H \in \Sn$, the matrix
\[
H_+ = (I-zz^T)H(I-zz^T) + zz^T
\]
satisfies the orthogonality condition $(H_+ - I) \perp (H_+ - H)$, and consequently
\[
\|(H-I)z\|^2 ~\le~ \|H-I\|^2 - \|H_+ - I\|^2,
\]
and furthermore $\lm(H_+) ~\ge~ \min \{ \lm(H),1 \}$.
\end{lem}

\pf
In Lemma \ref{orthogonal}, we consider the matrices $P = I-zz^T$ and $X = H-I$.  Then we have
\[
X_+ ~=~ (I-zz^T)(H-I)(I-zz^T) ~=~ H_+ - I,
\]
so the orthogonality condition follows.  Using this, and noting $H_+ z = z$, since
\[
\|(H-I)z\|^2 ~=~ \|(H-H_+)z\|^2 ~\le~ \|H- H_+\|^2 ~=~ \|H-I\|^2 - \|H_+ - I\|^2,
\]
and the first inequality follows.

Turning to the second inequality, choose a unit vector $u \in \Rn$ satisfying $u^T H_+ u = \lm(H_+)$.  Then we deduce
\begin{eqnarray*}
\lm(H_+)
&=&
\big( u - (z^T u)z \big)^T H \big( u - (z^T u)z \big) + (z^T u)^2 \\
& \ge &
\|u - (z^T u)z\|^2 \lm(H) + (z^T u)^2 \\
& = &
\big( 1 - (z^Tu)^2 \big) \lm(H) + (z^T u)^2  \\
& = &
\big( 1 - \lm(H) \big) (z^Tu)^2  + \lm(H).
\end{eqnarray*}
The result now follows.
\finpf

We can now complete our analysis in the quadratic case.

\begin{thm}
BFGS updating (Algorithm \ref{descent}) terminates for any strictly convex quadratic function 
$f$ at any noncritical point.
\end{thm}

\pf
We apply Lemma \ref{tool1}.  As we have argued, it suffices to consider Algorithm \ref{bfgs-norm} at any nonzero point $x \in \Rn$, so suppose by way of contradiction that the procedure does not terminate.  As the iterations progress, the nonnegative quantity $\|H-I\|$ is nonincreasing.  The first important consequence is the uniform boundedness of the matrix $H$ and hence that of the vector $s$.  Secondly, we also deduce $(H-I)z \to 0$.  Denoting the initial matrix $H$ by $H_0$, we see by induction the inequality 
$\lm(H) \ge \min \{ \lm(H_0),1 \} > 0$.
Thus the matrix $H^{-1}$ also stays uniformly bounded, so we deduce 
$\frac{1}{\|s\|}{s + x} = (I - H^{-1})z \to 0$.
Consequently we see $s \to -x$, contradicting the assumption that the procedure does not terminate.
\finpf

The argument above in fact proves more.

\begin{thm}
The linesearch-free BFGS method converges to the minimizer of any strictly convex quadratic function.
\end{thm}

\pf
As above, after a suitable change of variables it suffices to prove the result for the case 
$f(x) = \frac{1}{2}\|x\|^2$.  Starting from any nonzero initial point $x \in \Rn$ and matrix $H \in \Sn_{++}$, we are therefore repeating the following procedure.

{\rm
\begin{algorithmic}
\STATE
\WHILE{$\|x - Hx\| \ge \|x\|$}
\STATE	$s = -Hx$; $x_+ = x+s$; $z = \frac{s}{\|s\|}$; $H = (I-zz^T)H(I-zz^T) + zz^T$;
\ENDWHILE
\STATE	$x=x_+$;
\end{algorithmic}
}

\bigskip

\noindent
Exactly as before, we argue that the matrices $H$ and $H^{-1}$ remain uniformly bounded.  By the definition of the sequence of iterates $x$, we know $\|x\|$ is nonincreasing, so the steps $s$ are also uniformly bounded.  We deduce $x_+ \to 0$, since as before we know
\bmye \label{convergence}
\frac{s + x}{\|s\|}  \to 0.
\emye
If the iterates $x$ do not converge to zero, then they are uniformly bounded away from zero, and hence eventually we always accept the step because $\|x_+\| < \|x\|$.  But this is a contradiction, since 
$x_+ \to 0$.
\finpf

The argument shows a little more.  After taking the inner product with the unit vector $\frac{s}{\|s\|}$, equation (\ref{convergence}) shows $\frac{1}{\|s\|^2} s^T x \to -1$
and hence $\frac{1}{\|s\|^2} (\|x_+\|^2 - \|x\|^2) \to -1$.
Thus eventually we always accept the step because $\|x_+\| < \|x\|$.  We have thus shown that the linesearch-free BFGS method applied to a strictly convex quadratic always accepts the step eventually.  The method then reduces to the classical method, and hence converges superlinearly to the minimizer \cite{nocedal_wright}.

\section{BFGS updating for nonsmooth functions} \label{nonsmooth-updating}
In practice we can apply the classical BFGS method directly to nonsmooth functions, after a randomized initialization, as we noted in the introduction (see \cite{BFGS}).  However, our aim here is to illuminate the effect of the BFGS update for nonsmooth functions as simply as possible, so we first consider more formally how we should define it.

To that end, consider a convex function $f \colon \Rn \to \R$, possibly nonsmooth.  Given a current point $x$, subgradient $g \in \partial f(x)$, and matrix $H \in \Sn_{++}$, generalizing the classical BFGS method (with a unit step) leads to the following update:
\begin{myeqnarray}
s = -Hg, & &  x_+ = x+s \label{newupdate1} \\
g_+ &\in& \mbox{argmax} \big\{ z^Ts : z \in \partial f(x_+) \big\} \label{newupdate2} \\
y = g_+ - g, ~~~ V = I - \frac{sy^T}{s^Ty}, & & H_+ = VHV^T + \frac{ss^T}{s^Ty}.  \label{newupdate3}
\end{myeqnarray}
A priori it seems that we might choose the subgradient $g_+$ arbitrarily from $\partial f(x_+)$.   
The motivation for the particular choice in equation (\ref{newupdate2}) deserves some explanation.

As we have discussed, practical BFGS methods, both in the classical smooth case and in the nonsmooth case, use a line search, suitably scaling the step $s$ before updating $x \leftarrow x_+$,~ $g \leftarrow g_+$, $H \leftarrow H_+$ and repeating.  The curvature condition in the line search depends crucially on the directional derivative $f'(x_+;s)$ of the objective $f$ at the new point $x_+$ along the search direction $s$.  By standard convex analysis \cite{cov_lift}, that directional derivative is given by
\[
f'(x_+;s) ~=~ \max \big\{ z^Ts : z \in \partial f(x_+) \big\} ~=~ g_+^Ts.
\]
Thus our choice of the new subgradient $g_+$ corresponds to the correct linear approximation to the objective $f$ at the new point $x_+$ along the direction of the last step $s$.  Worth noting too is that, for analogous reasons, this choice of subgradient is also reminiscent of the subgradients generated by bundle methods for convex optimization \cite{hu-lem}.

Repeating this generalized unit-step BFGS update at a fixed point $x$ presents a fresh difficulty:  not only must we choose a subgradient $g_+ \in \partial f(x_+)$ but we may also update the original subgradient $g \in \partial f(x)$.  An analogous argument to the previous paragraph suggests the choice
\bmye \label{newupdate4}
g_{++} ~\in~  \mbox{argmax} \big\{ z^Ts : z \in \partial f(x) \big\},
\emye
since this corresponds to the correct linear approximation to the objective $f$ at the fixed point $x$ along the direction of the last trial step $s$:  $f'(x;s) = g_{++}^Ts$.  We are therefore led to the following generalized definition.

\begin{defn} \label{new-unit-step}
{\rm 
Consider a convex function $f \colon \Rn \to \R$ and a point $x$ in $\Rn$.  The {\em (nonsmooth) unit-step BFGS update} is the set-valued mapping 
\[
\mbox{BFGS}_{f,x} \colon \partial f(x) \times \Sn_{++} ~\tto~ \partial f(x) \times \Sn_{++}
\]
defined, for any subgradient $g \in \partial f(x)$ and a matrix $H \in \Sn_{++}$, by
\[
\mbox{BFGS}_{f,x}(g,H) ~=~ \big\{ (g_{++},H_+) : \mbox{(\ref{newupdate1}), (\ref{newupdate2}), (\ref{newupdate3}), (\ref{newupdate4})~ hold} \big\}.
\]
}
\end{defn}

\noindent
Notice that the set of updates $\mbox{BFGS}_{f,x}(g,H)$ is empty if $s^Ty = 0$ in equation (\ref{newupdate3}).  When the function $f$ is smooth, the set $\mbox{BFGS}_{f,x}(\nabla f(x),H)$ consists of just one element, namely the matrix we called $\mbox{BFGS}_{f,x}(H)$ in our previous notation.

We can now pose our questions at the end of the introduction more precisely.  In particular, consider 
the nonsmooth unit-step BFGS update algorithm for the objective function $f$ at a fixed point $x$:
\bmye \label{newdescent}
\mbox{\bf while}~ f(x - Hg) \ge f(x),~~  (g,H) \leftarrow \mbox{BFGS}_{f,x}(g,H).
\emye
\begin{itemize}
\item
If $x$ is not a minimizer, what conditions guarantee termination with descent:  
$f(x-Hg) < f(x)$?
\item
If $x$ is a minimizer, what conditions guarantee that the step $-Hg$ converges to zero?
\end{itemize} 

\section{BFGS updating for sublinear functions}
We now return to the interesting special case we discussed in Section \ref{space-dilation}, when the point of interest is $x=0$, and the function $f$ is sublinear.  In that case the unit-step BFGS update simplifies.  As a consequence of the following result, whose proof is an easy exercise, the distinction between the sets of acceptable subgradients 
$g_+$ and $g_{++}$ vanishes in this case.  To simplify the update in this case, we can choose 
$g_{++} = g_+$.

\begin{prop}
For any sublinear function $f \colon \Rn \to \R$ and any vector $s \in \Rn$, the maximum value of the linear function $\ip{s}{\cdot}$ over the subdifferential $\partial f(0)$ is $f(s)$, and the set of maximizers is $\partial f(s)$.
\end{prop}

We can consider any sublinear function $f \colon \Rn \to \R$ as the support function $\delta_C^*$ of a nonempty compact set $C \subset \Rn$, namely $C = \partial f(0)$. Hence the nonsmooth unit-step BFGS update algorithm (\ref{newdescent}) for deciding whether or not the point zero minimizes a sublinear function $f = \delta_C^*$ is equivalent to the following algorithm for deciding whether or not zero lies in the compact convex set $C$.
\newpage

\begin{alg}[BFGS for $0 \in C$]
\label{membership}
{\rm
\begin{algorithmic}
\STATE
\STATE  Choose $g \in C$, and $H \in \Sn_{++}$;
\FOR{$k=0,1,2,\dotsc$}
\IF{$g=0$}
\STATE  terminate with ``$0 \in C$'';
\ENDIF
\STATE  $s = -Hg$;
\STATE  Find a maximizer $g_+$ of $\ip{\cdot}{s}$ over $C$;
\IF{$g_+^Ts < 0$}
\STATE  terminate with $s$ ``normal to hyperplane separating $0$ from $C$'';
\ENDIF
\STATE  $y = g_+ - g$; $V = I - \frac{sy^T}{s^Ty}$; $H_+ = VHV^T + \frac{ss^T}{s^Ty}$; $H = H_+$; 
$g = g_+$;
\ENDFOR
\end{algorithmic}
}
\end{alg}

\noindent
Notice that if both the stopping conditions fail, so $g \ne 0$ and $g_+^T s \ge 0$, then
\[
y^T s = g_+^T s - g^Ts \ge g^T H g > 0,
\]
so the BFGS update is well-defined.

We can translate the questions at the end of the previous section for this special case.  Consider Algorithm \ref{membership} applied to a compact convex set $C$.  What conditions ensure the following properties?
\begin{itemize}
\item
$0 \not\in C ~\Rightarrow$ correct termination.
\item
$0 \in C ~\Rightarrow$ either correct termination or convergence of the step $s$ to zero.
\end{itemize}

The algorithm depends on being able to maximize linear functionals over the compact convex set 
$C \subset \Rn$, so is most realistic when $C$ is the convex hull of a possibly simpler (even finite) compact set $D \subset \Rn$.  In that case we can choose to restrict our attention to maximizers 
$g_+$ in $D$ rather than $C$.  Furthermore, if $0 \not\in D$, we can omit the first termination criterion.   We then arrive at the following algorithm for deciding whether or not zero lies in the convex hull of a compact set $D \subset \Rn$ not containing zero.
\newpage

\begin{alg}[BFGS for $0 \in \mbox{conv}\,D$]
\label{membership2}
{\rm
\begin{algorithmic}
\STATE
\STATE  Choose $g \in D$, and $H \in \Sn_{++}$;
\FOR{$k=0,1,2,\dotsc$}
\STATE  $s = -Hg$;
\STATE  Find a maximizer $g_+$ of $\ip{\cdot}{s}$ over $D$;
\IF{$g_+^Ts < 0$}
\STATE  terminate with $s$ ``normal to hyperplane separating $0$ from $\mbox{conv}\,D$'';
\ENDIF
\STATE  $y = g_+ - g$; $V = I - \frac{sy^T}{s^Ty}$; $H = VHV^T + \frac{ss^T}{s^Ty}$; $g = g_+$;
\ENDFOR
\end{algorithmic}
}
\end{alg}

Many authors (such as \cite[p.\ 1051]{bland-goldfarb-todd}) have noted the similarities between quasi-Newton algorithms like the BFGS method, and the Ellipsoid Algorithm and related space-dilation techniques (especially the Shor r-algorithm \cite[Section 3.6]{shor}).  The Ellipsoid Algorithm for minimizing, over the unit ball, the support function 
$\delta_C^*$, when the set $C$ is the convex hull of a compact set $D \subset \Rn$ not containing zero, takes the following form \cite[p.\ 249]{bubeck}.  We note the similarities with Algorithm~\ref{membership2}.

\begin{alg}[Ellipsoid algorithm for $0 \in \mbox{conv}\,D$]
{\rm
\begin{algorithmic}
\STATE
\STATE  $x = 0$; $H = I$;
\FOR{$k=0,1,2,\dotsc$}
\IF{$\|x\| > 1$}
\STATE  $g = x$;
\ELSE
\STATE  Find a maximizer $g$ of $\ip{\cdot}{x}$ over $D$;
\IF{$g^Tx < 0$}
\STATE  terminate with ``$x$ separates $0$ from $\mbox{conv}\,D$'';
\ENDIF
\ENDIF
\STATE  $s = -Hg$; $x = x + \frac{s}{(n+1)\sqrt{-s^Tg}}$;
$H = \frac{n^2}{n^2-1}\big(H + \frac{2ss^T}{(n+1)s^Tg}\big)$;
\ENDFOR
\end{algorithmic}
}
\end{alg}

\section{Symmetry and the unit ball}
We begin with the second of our two questions: how does Algorithm \ref{membership} behave when the compact convex set $C$ contains zero?  We recall the following measure of the symmetry of the set $C$:
\[
\mbox{sym}(C) ~=~ \max \{t : g \in C ~\Rightarrow -tg \in C \}.
\]
This measure often appears in complexity analysis for convex optimization \cite{nesterov-nemirovskii,renegar-book,freund-complexity}.  We also use the following standard result, whose proof entails simple linear algebra \cite[equation (6.45)]{nocedal_wright}.

\begin{lem}
\label{det}
The matrices $H$ and $H_+$ in Algorithm \ref{membership} satisfy
\[
\frac{\det H_+}{\det H} ~=~ -\frac{s^Tg}{s^Ty}.
\]
\end{lem}

\noindent
In fact this result holds for any matrices $H,H_+ \in \Sn$ related via the BFGS update equations (\ref{newupdate1}) and (\ref{newupdate3}).

Our next result shows that when the set $C$ contains zero in its interior, the determinant of the matrix $H$ must converge to zero, and at a linear rate controlled by the symmetry measure.

\begin{prop}
\label{smaller}
If the compact convex set $C$ contains zero, then the matrices $H$ and $H_+$ in Algorithm \ref{membership} always satisfy
\[
\det H_+ ~\le~ \frac{\det H}{1+\mbox{\rm sym}(C)}.
\]
\end{prop}

\pf
Since $g \in C$, by definition we have $-\mbox{sym}(C)g \in C$.  By Lemma \ref{det} we deduce
\[
\frac{\det H}{\det H_+} ~=~ \frac{s^T(g - g_+)}{s^Tg}
 ~=~ 1 + \frac{\max_C \ip{\cdot}{s}}{-s^Tg} ~\ge~ 1 + \frac{\ip{-\mbox{sym}(C)g}{s}}{-s^Tg}.
\]
The result follows.
\finpf

When the set $C$ is simply the unit ball, Algorithm \ref{membership} becomes particularly simple.  Numerical experiments suggest the following conjecture.

\begin{con}[BFGS for the unit ball]
\label{unit}
Given any initial unit vector $g \in \Rn$ and matrix $H \in \Sn_{++}$, if we repeatedly set
\[
s=-Hg,~~
g_+ = \frac{s}{\|s\|},~~
y = g_+ - g,~~
V = I - \frac{sy^T}{s^Ty},~~
H_+ = VHV^T + \frac{ss^T}{s^Ty},
\]
and update $g = g_+$ and $H = H_+$, then the trial step $s$ converges to zero.
\end{con}

\begin{figure}  
\begin{center}
\includegraphics[width=10cm]{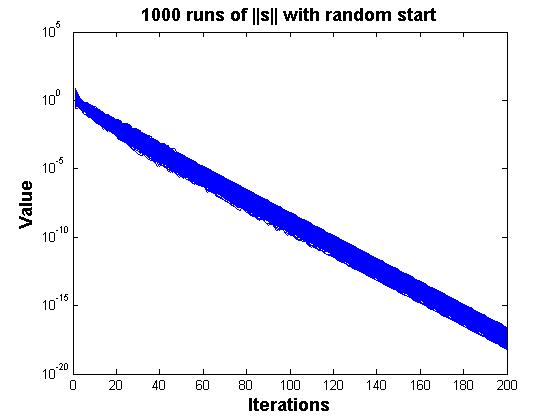}
\end{center}
\caption{1000 random runs of the iteration in Conjecture \ref{unit}.  Step size $\|s\|$ plotted against iteration count.}
\label{conjecture-evidence}
\end{figure}

\begin{figure}  
\begin{center}
\includegraphics[width=10cm]{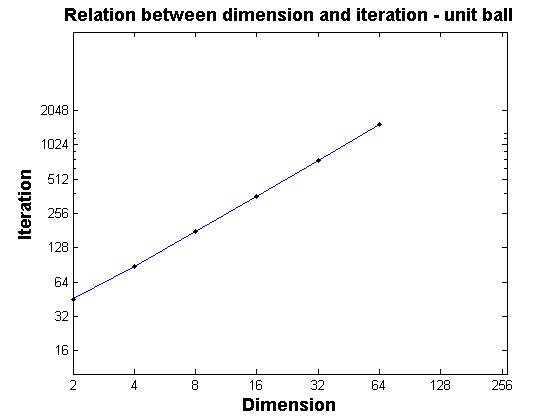}
\end{center}
\caption{Mean number of iterations, over 200 random runs, to reduce $\|s\|$ by a factor $10^{-8}$ in Conjecture \ref{unit}, plotted against dimension $n$.}
\label{dimension-dependence}
\end{figure}

\noindent
Figure \ref{conjecture-evidence} shows overlaid plots of $\|s\|$ against iteration count for a thousand randomly initiated runs in dimension $n=5$.  Such numerical results strongly suggest a linear convergence rate, and one that grows quite slowly with dimension $n$.  Figure \ref{dimension-dependence} plots against dimension $n$, on a log-log scale, the number of iterations (averaged over 200 random runs) to reduce $\|s\|$ by a factor $10^{-8}$ in Conjecture \ref{unit}:  the number grows roughly like $n^{1/\sqrt{2}}$.

As a first theoretical step we prove the following result.  We denote the largest and smallest eigenvalues of $H$ by 
$\lambda_{\max}(H)$ and $\lambda_{\min}(H)$ respectively, and we write $E \succ F$ for matrices $E,F \in \Sn$ to mean $E-F \in \Sn_{++}$. 

\begin{thm}
The matrices $H$ and $H_+$ in Conjecture \ref{unit} satisfy
\[
\det H_+ \le \frac{1}{2} \det H ~~\mbox{and}~~  \lambda_{\max}(H_+) \le \lambda_{\max}(H).
\]
\end{thm}

\pf
The first inequality follows from Proposition \ref{smaller}.  Turning to the second, notice that the update of the matrix $H$ to $H_+$ is positively homogeneous:  if we replace $H$ by $\gamma H$ for some positive scalar $\gamma$ then $H_+$ is replaced by $\gamma H_+$.  After a suitable scaling of $H$, we can therefore assume that the step $s$ is a unit vector.  Similarly, if we replace $H$ by $U^T H U$, for an orthogonal matrix $U$, and $g$ by $U^Tg$, then $H_+$ is replaced by $U^T H_+ U$.  After such an orthogonal transformation, we can therefore also assume that $s$ is just $e_1$, the first unit vector.  

Partitioning vectors, we can write $s = [1~0]^T$ and $g = -[\alpha~b]^T$,
for some scalar $\alpha \in (0,1]$ and vector $b \in \R^{n-1}$ satisfying $\alpha^2 + \|b\|^2 = 1$.  Since $s = -Hg$, we can also partition the matrix $H^{-1}$ as
\[
H^{-1} =
\left[
\begin{array}{cc}
\alpha & b^T \\
b & E 
\end{array}
\right],
\]
where, using the Schur complement, the matrix $E \in {\mathbf S}^{n-1}$ satisfies 
$E \succ \frac{bb^T}{\alpha}$.  We can write the BFGS update equivalently as
\[
H_+^{-1} = H^{-1} + \frac{yy^T}{s^Ty} + \frac{gg^T}{s^Tg}
\]
(see \cite{nocedal_wright}), and we deduce
\[
H_+^{-1} = 
\left[
\begin{array}{cc}
1+ \alpha & b^T \\
b & E - \frac{bb^T}{\alpha(1+\alpha)} 
\end{array}
\right].
\]

For any positive scalar $\lambda$ satisfying $\lambda < \lambda_{\min}(H^{-1}) = \big(\lambda_{\max}(H)\big)^{-1}$, we seek to show  
$\lambda < \lambda_{\min}(H_+^{-1})  = \big(\lambda_{\max}(H_+)\big)^{-1}$.  Equivalently, via the Schur complement, we know
\[
\alpha - \lambda > 0 ~~\mbox{and}~~  E - \lambda I \succ \frac{bb^T}{\alpha - \lambda},
\]
and we seek to prove
\[
1 + \alpha - \lambda > 0 ~~\mbox{and}~~  
E  - \frac{bb^T}{\alpha(1+\alpha)} - \lambda I \succ \frac{bb^T}{1 + \alpha - \lambda}.
\]
The first inequality is immediate.  The second follows from the inequality
\[
\frac{1}{\alpha - \lambda} \ge \frac{1}{\alpha(1+\alpha)} + \frac{1}{1 + \alpha - \lambda},
\]
which is equivalent to the inequality 
$\alpha(1+\alpha) \ge (\alpha - \lambda)(1 + \alpha - \lambda)$, a consequence of the monotonicity of the left-hand side.
\finpf

\section{Cholesky factors and line segments} \label{cholesky-factors}
The Shor r-algorithm (see equation (\ref{shor-factors})) uses the search direction $s = -V^TV g$, where $g$ is a current subgradient.  The analogue of the quasi-Newton matrix $H$ is $V^TV$, and the method updates the factor $V$.  As we commented at the end of Section~\ref{space-dilation}, we can take a similar approach to the BFGS algorithm.  Instead of updating the inverse Hessian approximation $H$ directly through the BFGS formula (\ref{update2}), it can be useful (see \cite{powell87,byatt-coope-price}) to update a factored form $H = T^T T$, where the matrix $T$ is invertible.  In that case we can write the update as 
$H_+ = T_+^T T_+$, where
\[
T_+ = T(I - qs^T),~~  \mbox{for}~~ q = \frac{y}{s^Ty} + \frac{g}{\sqrt{-s^Tg s^Ty}}.
\]

Consider the BFGS algorithm for $0 \in C$, with this notation.  After some algebra and the change of variables $h = Tg$, $p = Tg_+$ (where $g_+$ is the updated vector $g$), and 
$P = TC$, we arrive at the following 
algorithm for deciding whether or not zero lies in a compact convex set $P$.

\begin{alg}[Cholesky BFGS for $0 \in P$]
{\rm
\begin{algorithmic}
\STATE
\STATE  Choose $h \in P$;
\FOR{$k=0,1,2,\dotsc$}
\IF{$h=0$}
\STATE  terminate with ``$0 \in P$'';
\ENDIF
\STATE  Find a minimizer $p$ of $\ip{\cdot}{h}$ over $P$;
\IF{$p^Th > 0$}
\STATE  terminate with ``$0 \not\in P$'';
\ENDIF
\STATE  $e = h-p$; $\beta = h^Te$; $W = I - \frac{eh^T}{\beta} + \frac{hh^T}{\|h\|\sqrt{\beta}}$;
$P = WP$; $h = Wp$;
\ENDFOR
\end{algorithmic}
}
\end{alg}

As an example, consider the convex hull of nonzero vectors $a_i \in \Rn$ indexed by a finite set 
$I$.  We can implement the algorithm above as follows.

\begin{alg}[Cholesky BFGS for $0 \in \mbox{conv}\{a_i : i \in I\}$]
\label{cholesky}
{\rm
\begin{algorithmic}
\STATE
\STATE  Choose $i \in I$;
\FOR{$k=0,1,2,\dotsc$}
\STATE  Find $j \in I$ minimizing $a_i^T a_j$;
\IF{$a_i^T a_j > 0$}
\STATE  terminate with ``$0$ lies outside the convex hull'';
\ENDIF
\STATE  $e = a_i - a_j$; $\beta = a_i^Te$;
\FOR{each $r \in I$}
\STATE  $a_r = a_r - (a_i^Ta_r)(\frac{e}{\beta} - \frac{a_i}{\|a_i\|\sqrt{\beta}})$;
\ENDFOR
\STATE  $i = j$;
\ENDFOR
\end{algorithmic}
}
\end{alg}

To illustrate, consider how this method behaves for a set of just two distinct nonzero vectors $a_1 = c \ne d = a_2$ in $\Rn$.  The algorithm becomes the following.

\begin{alg}[Cholesky BFGS for $0 \in \mbox{$[c,d]$}$]
\label{segment}
{\rm
\begin{algorithmic}
\STATE
\FOR{$k=0,1,2,\dotsc$}
\IF{$c^Td > 0$}
\STATE  terminate with ``$0 \notin [c,d]$'';
\ENDIF
\STATE  $e = c - d$; $\beta = c^Te$;
\STATE  $d_+ = d - (c^Td)(\frac{e}{\beta} - \frac{c}{\|c\|\sqrt{\beta}})$;
$c_+ = c - (c^Tc)(\frac{e}{\beta} - \frac{c}{\|c\|\sqrt{\beta}})$;
\STATE  $c=d_+$; $d=c_+$;
\ENDFOR
\end{algorithmic}
}
\end{alg}

We introduce a measure to track the conditioning of the line segments:
\[
\gamma[c,d] ~=~ \frac{\sqrt{\|c\|^2\|d\|^2 - (c^Td)^2}}{\|c-d\|^2}.
\]
This quantity is well-defined since the right-hand side is symmetric in $c$ and $d$.  It is also invariant under scaling and orthogonal transformations: 
$\gamma(\alpha[c,d]) = \gamma[c,d] = \gamma(U[c,d])$ for any nonzero scalar $\alpha$ and any $n$-by-$n$ orthogonal matrix $U$.

Obviously the vectors in the algorithm all evolve in the two-dimensional space spanned by the original line segment $[c,d]$.  Choosing a suitable basis, we therefore lose no generality in studying the special case $a_1 = c = [1~0]^T$ and $a_2 = d = [-p~q]^T$ with $q \ge 0$,
in which case we have
\bmye \label{gamma}
\gamma[c,d] = \frac{q}{(1+p)^2 + q^2}.
\emye
We arrive at the following tool for recognizing when the algorithm will terminate.

\begin{lem}[Angle recognition]
Assuming the condition $c^Td \le 0$, the line segment $[c,d]$ satisfies $\gamma[c,d] \le \frac{1}{2}$.  Under the additional assumption $0 \not\in [c,d]$, the segment also satisfies $\gamma[c,d] > 0$.
\end{lem}

\pf
Assuming the special case above, we have $p \ge 0$.  We deduce
\[
(1+p)^2 + (q-1)^2 \ge 1,
\]
and the first claim now follows from equation (\ref{gamma}).  The second claim is easy.
\finpf

If the original line segment $[c,d]$ does not contain zero, eventually the termination criterion will hold, as a consequence of the following conditioning improvement.

\begin{lem}
If $c^Td \le 0$, then $\gamma[c_+,d_+] \ge \gamma[c,d] + (\gamma[c,d])^3$.
\end{lem}

\pf
We consider the special case above again, so by assumption, $p \ge 0$.  Since (\ref{gamma}) holds, in particular we have
\bmye  \label{upper}
\gamma[c,d] \le \frac{q}{(1+p)^{3/4}}.
\emye
A quick calculation shows 
\[
\gamma[c_+,d_+] = \frac{q}{(1+p)^{3/2}}.
\]
We deduce
\[
\frac{\gamma[c_+,d_+]}{\gamma[c,d]} = 
\frac{(1+p)^2 + q^2}{(1+p)^{3/2}} =
(1+p)^{1/2} + \frac{q^2}{(1+p)^{3/2}} \ge
1 + (\gamma[c,d])^2,
\]
by inequality (\ref{upper}).  The result follows.
\finpf

We can be more precise, using the following lemma.

\begin{lem}
For any $K > 0$, consider the finite sequence $(\beta_k)$ defined (for integers $k \ge 0$) by 
$\beta_k =(K+1-k)^{-1/2}$ for all $k \le K$.
Suppose a second sequence $(\gamma_k)$ satisfies $\gamma_0 \ge \beta_0$ and 
$\gamma_{k+1} \ge \gamma_k + \gamma_k^3$ for all $k \le K-1$.
Then $\gamma_k \ge \beta_k$ for all $k \le K$.
\end{lem}

\pf
To prove the result by induction, we just need to show
$\beta_{k+1} \le \beta_k + \beta_k^3$ for all $k \le K-1$.
Squaring both sides, we obtain 
\[
\beta_{k+1}^2 = \frac{\beta_k^2}{1-\beta_k^2} \le \beta_k^2 + 2\beta_k^4 + \beta_k^6,
\]
or equivalently,
$\beta_k^4 + \beta_k^2 \le 1$.  This last inequality is valid, since $\beta_k^2 \le \frac{1}{2}$.
\finpf

Suppose the original segment $[c,d]$ does not contain zero, and denote the condition measure 
$\gamma[c,d]$ after $k = 0,1,2,\ldots$ iterations by $\gamma_k$.  Then we have
$\gamma_0 = \gamma[c,d]$ and $\gamma_{k+1} \ge \gamma_k + \gamma_k^3$
so we deduce by the previous lemma,
$\gamma_k \ge (\gamma_0^{-2}-k)^{-1/2}$ for all $k \le \gamma_0^{-2} - 1$.
In particular, this inequality holds when $k$ is the integral part 
$\lfloor \gamma_0^{-2} - 1 \rfloor$, if the algorithm has not already terminated.  In that case,
$k >  \gamma_0^{-2} - 2$, so 
$
\gamma_k \ge (\gamma_0^{-2}-k)^{-1/2} > 2^{-1/2} > \frac{1}{2}
$,
so the algorithm terminates.

We have proved the following result.

\begin{thm}
For any distinct nonzero vectors $c,d \in \Rn$, if the line segment $[c,d]$ does not contain zero, then, after a number of iterations not exceeding 
\[
\frac{\|c-d\|^4}{\|c\|^2\|d\|^2 - (c^Td)^2}
\]
Algorithm \ref{segment} terminates correctly.
\end{thm}

\section{Conclusion}
This work explores the relative effectiveness of the BFGS method and the Shor r-algorithm in the context of nonsmooth convex optimization. Incorporating line searches complicates the analysis, so here we try to separate their impact from the effect of the quasi-Newton or Shor update.  In particular, we consider a  simple linesearch-free BFGS algorithm.

Our experiments illustrate the effectiveness of improving the local metric in nonsmooth optimization.  We focus especially on simple examples where the current subdifferential is a polytope, ball or ellipsoid, presenting both numerical and theoretical results.  The algorithms simplify even further, conceptually, when rather than updating the approximate Hessian matrix, we instead work with its Cholesky factors.  In summary, this exploration only heightens our appreciation for the mysterious power of the BFGS methodology.

\bibliographystyle{plain}
\small
\parsep 0pt

\def\cprime{$'$} \def\cprime{$'$}


\end{document}